\documentclass[11pt,reqno]{amsart}
\UseRawInputEncoding

\usepackage{latexsym,bm, amsfonts,amsthm,amsmath,mathrsfs,CJK}

\makeatletter

\newcommand{\Rmnum}[1]{\expandafter\@slowromancap\romannumeral #1@}
\makeatother

 \newtheorem{lem}{Lemma}[section]  \newtheorem{thm}{Theorem}[section]
 \newtheorem{defn}{Definition}[section] \newtheorem{rmk}{Remark}[section] 
\numberwithin{equation}{section}

\newcommand{\dif}{\mathrm{d}} \DeclareMathAlphabet{\mathsfsl}{OT1}{cmss}{m}{sl} \DeclareMathAlphabet{\mathpzc}{OT1}{pzc}{m}{it}

    \newcommand{\ee}{\mathbb{E}}

\newcommand{\nn}{\mathbb{N}}
\newcommand{\rr}{\mathbb{R}}
    
\newcommand{\vv}{\mathbb{V}}
 \def\CC{\mathcal C}   \def\FF{\mathcal F}  \def\HH{\mathcal H}

\def\d"{^{\prime\prime}} \def\d'{^{\prime}}

\begin{document}

		\title[]{Note on the complete moment convergence for moving average process of a class of random variables under sub-linear expectations}\thanks{Project supported by Science and Technology Research Project of Jiangxi Provincial Department of Education of China (Nos. GJJ2201041), Doctoral Scientific Research Starting Foundation of Jingdezhen Ceramic University ( Nos.102/01003002031).}
		\date{} \maketitle


		
		\begin{center}
			XU Mingzhou~\footnote{Email: mingzhouxu@whu.edu.cn}
			\\
			School of Information Engineering, Jingdezhen Ceramic University\\
			Jingdezhen 333403, China
		\end{center}

		\renewcommand{\abstractname}{~}
		
		{\bf Abstract:}
		In this paper, the complete moment convergence for the partial sums of moving average processes $\{X_n=\sum_{i=-\infty}^{\infty}a_iY_{i+n},n\ge 1\}$ is proved under some proper conditions, where $\{Y_i,-\infty<i<\infty\}$ is a doubly sequence of identically distributed, negatively dependent random variables under sub-linear expectations and $\{a_i,-\infty<i<\infty\}$ is an absolutely summable sequence of real numbers. The results established in sub-linear expectation spaces generalize the corresponding ones in probability space.
		
		
		\vspace{-3mm}

		\section{Introduction }
		
		In order to study the uncertainty in probability, Peng \cite{peng2007g,peng2019nonlinear} initiated  the concepts of the sub-linear expectations space.Inspired by the works of Peng \cite{peng2007g,peng2019nonlinear}, many scholars try to prove the results under sub-linear expectations space, extending the corresponding ones in classic probability space. Zhang\cite{zhang2015donsker,zhang2016exponential,zhang2016rosenthal} established Donsker's invariance principle, exponential inequalities and Rosenthal's inequality under sub-linear expectations. Under sub-linear expectations, Xu et al. \cite{xu2022complete}, Xu and Kong \cite{xu2023noteon} studied complete convergence for weighted sums of negatively dependent random variables. For more limit theorems under sub-linear expectations, the interested readers could refer to Xu and Zhang \cite{xu2019three,xu2020law}, Wu and Jiang \cite{wu2018strong},  Zhang and Lin \cite{zhang2018marcinkiewicz}, Zhong and Wu \cite{zhong2017complete}, Chen \cite{chen2016strong}, Chen and Wu \cite{chen2022complete}, Zhang \cite{zhang2022strong}, Hu, Chen and Zhang \cite{hu2014big}, Gao and Xu \cite{gao2011large}, Kuczmaszewska \cite{kuczmaszewska2020complete}, Xu and Cheng \cite{xu2021convergence,xu2022small,xu2022note}, Xu et al. \cite{xu2022complete}, Yu and Wu \cite{yu2018Marcinkiewicz} and references therein.
		
		In probability space, Hsu and Robbins \cite{hsu1947complete} introduced concept of complete convergence, Chow \cite{chow1988on}  studied complete moment convergence for independent random variables.  Li and Zhang \cite{li2004complete} got complete moment convergence of moving-average processes under dependence assumptions. There are vast of literatures on complete moment convergeces. For references on complete moment convergence in linear expectation space, the interested reader could refer to Baum and Katz \cite{baum1965convergence},  Sung \cite{sung2009moment}, Liu and Lin \cite{liu2006precise}, Zhou \cite{zhou2010complete},
		Qiu and Chen \cite{qiu2014acomplete,qiu2014bcomplete}, Liang et al. \cite{liang2010complete}, Shen et al. \cite{shen2016complete}, Guo and Zhu \cite{guo2013complete}, Kim and Ko \cite{kim2008complete}, Ko \cite{ko2015complete}, Meng et al. \cite{meng2021convergence},  Hosseini and Nezakati \cite{hosseini2019complete}, and refercences therein.  Recently Zhang and Ding \cite{zhang2017further} studied the complete moment convergence of the partial sums of moving average processes under some proper assumptions. Encouraged by the works of Zhang and Ding \cite{zhang2017further}, Chen and Wu \cite{chen2022complete}, and Xu et al. \cite{xu2022complete}, we try to prove the complete moment convergence for the partial sums of moving average processes generated by identically distributed, negatively dependent random variables under sub-linear expectations, extending the corresponding results of Zhang and Ding \cite{zhang2017further} in classic probability space. Our main contribution is that we extend the results of Zhang and Ding \cite{zhang2017further}, Li and Zhang \cite{li2004complete} in classic probability space to that under sub-linear expectations, and our method of proof here is heuristically inspired by that of Chen and Wu \cite{chen2022complete}, Zhang \cite{zhang2016rosenthal}, and Zhang and Ding \cite{zhang2017further}, which is different from that of Zhang and Ding \cite{zhang2017further}.
		
		We construct the remainders of this paper as follows. We present necessary basic notions, concepts and relevant  properties, and cite necessary lemma under sub-linear expectations in the next section. In Section 3, we present our main results, Theorem \ref{thm2.1}, Theorem \ref{thm2.2},  the proofs of which are given in Section 4.
		
		\section{Preliminary}
		\setcounter{equation}{0}
		We use similar notations as in the work by \cite{peng2019nonlinear}, Chen \cite{chen2016strong}, Zhang \cite{zhang2016rosenthal}. Assume that $(\Omega,\FF)$ is a given measurable space. Suppose that $\HH$ is a subset of all random variables on $(\Omega,\FF)$ such that  $X_1,\cdots,X_n\in \HH$ implies $\varphi(X_1,\cdots,X_n)\in \HH$ for each $\varphi\in \CC_{l,Lip}(\rr^n)$, where $\CC_{l,Lip}(\rr^n)$ stands for the linear space of (local lipschitz) function $\varphi$ satisfying
		$$
		|\varphi(\mathbf{x})-\varphi(\mathbf{y})|\le C(1+|\mathbf{x}|^m+|\mathbf{y}|^m)(|\mathbf{x}-\mathbf{y}|), \forall \mathbf{x},\mathbf{y}\in \rr^n
		$$
		for some $C>0$, $m\in \nn$ depending on $\varphi$.
		\begin{defn}\label{defn1} A sub-linear expectation $\ee$ on $\HH$ is a functional $\ee:\HH\mapsto \bar{\rr}:=[-\infty,\infty]$ fulfilling the following properties: for all $X,Y\in \HH$, we have
			\begin{description}
				\item[\rm (a)] If $X\ge Y$, then $\ee[X]\ge \ee[Y]$;
				\item[\rm (b)] $\ee[c]=c$, $\forall c\in\rr$;
				\item[\rm (c)] $\ee[\lambda X]=\lambda\ee[X]$, $\forall \lambda\ge 0$;
				\item[\rm (d)] $\ee[X+Y]\le \ee[X]+\ee[Y]$ whenever $\ee[X]+\ee[Y]$ is not of the form $\infty-\infty$ or $-\infty+\infty$.
			\end{description}
			
		\end{defn}
		A set function $V:\FF\mapsto[0,1]$ is said to be a capacity if
		\begin{description}
			\item[\rm (a)]$V(\emptyset)=0$, $V(\Omega)=1$;
			\item[\rm (b)]$V(A)\le V(B)$, $A\subset B$, $A,B\in \FF$.\\
		\end{description}
		A capacity $V$ is called sub-additive if $V(A+B)\le V(A)+V(B)$, $A,B\in \FF$.

		In this sequel, given a sub-linear expectation space $(\Omega, \HH, \ee)$, write $\vv(A):=\inf\{\ee[\xi]:I_A\le \xi, \xi\in \HH\}$, $\forall A\in \FF$ (see (2.3) and the definitions of $\vv$ above (2.3) in Zhang \cite{zhang2016exponential}. $\vv$ is a sub-additive capacity. Set
		$$
		C_{\vv}(X):=\int_{0}^{\infty}\vv(X>x)\dif x +\int_{-\infty}^{0}(\vv(X>x)-1)\dif x.
		$$
		As in 4.3 of Zhang \cite{zhang2016exponential}, hereafter, define an extension of $\ee$ on the space of all random variables by
		$$
		\ee^{*}(X)=\inf\left\{\ee[Y]:X\le Y,Y\in\HH\right\}.
		$$
		Then $\ee^{*}$ is a sublinear expectation on the space of all random variables, $\ee[X]=\ee^{*}[X]$, $\forall X\in \HH$, and $\vv(A)=\ee^{*}(I_A)$, $\forall A\in \FF$.

		Assume that $\mathbf{X}=(X_1,\cdots, X_m)$, $X_i\in\HH$ and $\mathbf{Y}=(Y_1,\cdots,Y_n)$, $Y_i\in \HH$  are two random vectors on  $(\Omega, \HH, \ee)$. $\mathbf{Y}$ is said to be negatively dependent to $\mathbf{X}$, if for each  $\psi_1\in \CC_{l,Lip}(\rr^m)$, $\psi_2\in  \CC_{l,Lip}(\rr^n)$, we have $\ee[\psi_1(\mathbf{X})\psi_2(\mathbf{Y})]\le\ee[\psi_1(\mathbf{X})] \ee[\psi_2(\mathbf{Y})]$ whenever $\psi_1(\mathbf{X})\ge 0$, $\ee[\psi_2(\mathbf{Y})]\ge 0 $, $\ee[\psi_1(\mathbf{X})\psi_2(\mathbf{Y})]<\infty$, $\ee[|\psi_1(\mathbf{X})|]<\infty$, $\ee[|\psi_2(\mathbf{Y})|]<\infty$, and either $\psi_1$ and $\psi_2$ are coordinatewise nondecreasing or $\psi_1$ and $\psi_2$ are coordinatewise nonincreasing (see Definition 2.3 of Zhang \cite{zhang2016exponential}, Definition 1.5 of Zhang \cite{zhang2016rosenthal}, Definition 2.5 in Chen \cite{chen2016strong}).
		$\{X_n\}_{n=1}^{\infty}$ is called to be a sequence of negatively dependent random variables, if $X_{n+1}$ is negatively dependent to $(X_1,\cdots,X_n)$ for each $n\ge 1$.
		
		Suppose that $\mathbf{X}_1$ and $\mathbf{X}_2$ are two $n$-dimensional random vectors defined, respectively, in sub-linear expectation spaces $(\Omega_1,\HH_1,\ee_1)$ and $(\Omega_2,\HH_2,\ee_2)$. They are called identically distributed if  for every $\psi\in \CC_{l,Lip}(\rr^n)$ such that $\psi(\mathbf{X}_1)\in \HH_1, \psi(\mathbf{X}_2)\in \HH_2$,
		$$
		\ee_1[\psi(\mathbf{X}_1)]=\ee_2[\psi(\mathbf{X}_2)], \mbox{  }
		$$
		whenever the sub-linear expectations are finite. $\{X_n\}_{n=1}^{\infty}$ is called to be identically distributed if for each $i\ge 1$, $X_i$ and $X_1$ are identically distributed.
		
		In the paper we suppose that $\ee$ is countably sub-additive, i.e., $\ee(X)\le \sum_{n=1}^{\infty}\ee(X_n)$, whenever $X\le \sum_{n=1}^{\infty}X_n$, $X,X_n\in \HH$, and $X\ge 0$, $X_n\ge 0$, $n=1,2,\ldots$. Hence $\ee^{*}$ is also countably sub-additive. Let $C$ stand for a positive constant which may change from place to place. $I(A)$ or $I_A$ represent the indicator function of $A$.
		\begin{defn}\label{defn01} A real valued function $l(x)$, positive and measurable on $[0,\infty)$, is said to be slowly varying at infinity if for each $\lambda>0$, $\lim_{x\rightarrow \infty}\frac{l(\lambda x)}{l(x)}=1$.
		\end{defn}
		
		We cite the following Rosenthal's inequality under sub-linear expectations.
		\begin{lem}\label{lem01}(See Theorem 2.1 of Zhang \cite{zhang2016rosenthal}) Let $\{X_n;n\ge 1\}$ be a sequence of negatively dependent random variables under sub-linear expectation space $(\Omega,\HH,\ee)$. Then there exists a positive constant $C=C_p$ depending on $p$ such that for $n\ge 1$, and $ p\ge 2$,
			\begin{eqnarray}\label{01}
				\nonumber&&	\ee\left[\max_{0\le k\le n}\left|\sum_{i=1}^{k}X_k\right|^p\right]\\
				&&\quad	\le C_p \left\{\sum_{k=1}^{n}\ee[|X_k|^p]+\left(\sum_{k=1}^{n}\ee[|X_k|^2]\right)^{p/2}+\left(\sum_{i=1}^{n}\left(-\ee(-X_i)\right)^{-}+\left(\ee(X_i)\right)^{+}\right)^p\right\},
			\end{eqnarray}
			where $a^{+}=\max\{a,0\}$, $a^{-}=\max\{-a,0\}$, $a\in \rr$.
		\end{lem}
		
		In the rest of this paper, assume that $\{Y_i,-\infty<i<\infty\}$ is a sequence of negatively dependent random variables, identically distributed as $Y$ with $\ee(Y)=\ee(-Y)=0$ under sub-linear expectation space $(\Omega,\HH,\ee)$, and $\{a_i,-\infty<i<\infty\}$ is a sequence of real numbers with $\sum_{i=-\infty}^{\infty}|a_i|<\infty$, and the moving average process $\{X_n,n\ge 1\}$ is defined by $X_n=\sum_{i=-\infty}^{\infty}a_iY_{i+n}$.
		
		\section{Main Results}
		\setcounter{equation}{0}
		Our main results are the following.
		\begin{thm}\label{thm2.1} Assume that $l$ is a function slowly varying at infinity. Let $\{a_i,-\infty<i<\infty\}$ be an absolutely summable sequence of real numbers. Suppose that $\{g(n);n\ge 1\}$ and $\{f(n);n\ge 1\}$ are two sequences of positive constants such that, for some $r\ge \max\{2,p\}$, $p\ge 1$,
			\begin{description}
				\item[\rm (C1)]  $f(n)\uparrow \infty$, $\frac{n}{f^p(n)}\rightarrow 0$;
				\item[\rm (C2)] $\sum_{m=1}^{k}\log\left(\frac{f(m+1)}{f(m)}\right)\sum_{n=1}^{m}\frac{ng(n)l(n)}{f(n)}=O(f^{p-1}(k)l(k))$;
				\item[\rm (C3)] $\sum_{m=k}^{\infty}\left[f^{1-r}(m)-f^{1-r}(m+1)\right]\sum_{n=1}^{m}\frac{ng(n)l(n)}{f(n)}=O(f^{p-r}(k)l(k))$;
				\item[\rm (C4)] $\sum_{m=1}^{k}\left[f(m+1)-f(m)\right]\sum_{n=1}^{m}\frac{ng(n)l(n)}{f(n)}=O(f^{p}(k)l(k))$;
				\item[\rm (C5)] $\sum_{m=1}^{\infty}\left[f^{1-r}(m)-f^{1-r}(m+1)\right]f^t(m+1)\sum_{n=1}^{m}\frac{n^{r/2}g(n)l(n)}{f(n)}<\infty$, where $t=\max\{0,2-p\}r/2$;
				\item[\rm (C6)] $\sum_{m=1}^{\infty}\left[f(m+1)-f(m)\right]f^{t'}(m+1)\sum_{n=1}^{m}\frac{n^{r/2}g(n)l(n)}{f(n)}<\infty$, where $t'=-\min\{2,p\}r/2$.
			\end{description}
			
			Suppose that $\{X_n=\sum_{i=-\infty}^{\infty}a_iY_{i+n},n\ge 1\}$ is a moving average process generated by a sequence of negatively dependent random variables $\{Y_i,-\infty<i<\infty\}$, identically distributed as $Y$ with $\ee(Y)=\ee(-Y)=0$ and $C_{\vv}\left(|Y|^p\left(1\vee l(f^{-1}(Y))\right)\right)<\infty$ under sub-linear expectation space $(\Omega,\HH,\ee)$, where $f^{-1}$ is the inverse function of $f$. Then, for all $\varepsilon>0$,
			\begin{equation}\label{2.1}
				\sum_{n=1}^{\infty}\frac{g(n)l(n)}{f(n)}C_{\vv}\left\{\left(\max_{1\le k\le n}\left|\sum_{j=1}^{k}X_j\right|-\varepsilon f(n)\right)^{+}\right\}<\infty.
			\end{equation}
		\end{thm}
		\begin{rmk}\label{rmk1}
			In Theorem \ref{thm2.1}, let $1\le s <2$, $p>s$, $f(n)=n^{1/s}$, $g(n)=n^{p/s-2}$, $r>\max\left\{2,p,\frac{p/s-1}{1/s-1/2},\frac{p/s-1}{p/2-1/2},\frac{p/s-1}{1/2}\right\}$, and assume that $\ee(Y)=\ee(-Y)=0$, $C_{\vv}\{|Y|^ph(|Y|^s)\}<\infty$. Then the conditions \rm(C1)-\rm(C6) mean that $n/n^{p/s}\rightarrow 0$, $\sum_{m=1}^{k}m^{-1}\sum_{n=1}^{m}n^{p/s-1-1/s}l(n)\approx Ck^{p/s-1/s}l(k)=O(k^{(p-1)/s}l(k))$, $\sum_{m=k}^{\infty}m^{(1-r)/s-1}\sum_{n=1}^{m}n^{p/s-1-1/s}l(n)\approx Ck^{p/s-r/s}l(k)=O(k^{(p-r)/s}l(k)) $,
			$\sum_{m=1}^{k}m^{1/s-1}\sum_{n=1}^{m}n^{p/s-1-1/s}l(n)\approx Ck^{p/s}l(k)=O(k^{p/s}l(k))$,\\ $\sum_{m=1}^{\infty}m^{(1-r)/s-1}(m+1)^{t/s}\sum_{n=1}^{m}n^{r/2+p/s-2-1/s}l(n)=\sum_{m=1}^{\infty}m^{p/s-2+t/s+r/2-r/s}l(m)<\infty$, $\sum_{m=1}^{\infty}m^{1/s-1}(m+1)^{t'}\sum_{n=1}^{m}n^{r/2+p/s-2-1/s}l(n)\approx \sum_{m=1}^{\infty}m^{-\min\{2,p\}r/2+p/s-2+r/2}l(m)<\infty$,
			and therefore we get
			$$
			\sum_{n=1}^{\infty}n^{p/s-2-1/s}l(n)C_{\vv}\left\{\left(\max_{1\le k\le n}|S_k|-\varepsilon n^{1/s}\right)^{+}\right\}<\infty,  \mbox{  for every $\varepsilon>0$,}
			$$
			which extends Theorem 1.1 of Li and Zhang \cite{li2004complete} in probability space to that under sub-linear expectations space.
			
		\end{rmk}
		\begin{rmk}\label{rmk01}
			In Theorem \ref{thm2.1}, let $1\le q<2$, $\alpha> 3/(q+2/q)>1$, $f(n)=n^{1/q}$, $g(n)=n^{\alpha-2}$, $l(n)=1$, $p=q\alpha$, and suppose that $r\ge \max\{2,q\alpha\}$ is large enough. Then the assumptions \rm(C1)-\rm(C6) imply respectively $n/n^{\alpha}\rightarrow 0$, $\sum_{m=1}^{k}m^{-1}\sum_{n=1}^{m}m^{\alpha -1-1/q}\approx Ck^{\alpha-1/q}=O(k^{\alpha-1/q})$, $\sum_{m=k}^{\infty}m^{(1-r)/q-1}\sum_{n=1}^{m}n^{\alpha-1-1/q}\approx Ck^{\alpha-r/q}=O(k^{\alpha-r/q})$, \\
			$\sum_{m=1}^{k}m^{1/q-1}\sum_{n=1}^{m}n^{\alpha-1-1/q}\approx Ck^{\alpha}=O(k^{\alpha})$,
			\begin{eqnarray*} &&\sum_{m=1}^{\infty}m^{(1-r)/q-1}m^{\max\{0,2-q\alpha\}r/2}\sum_{n=1}^{m}n^{r/2+\alpha-2-1/q},\\
				&&\quad\le \begin{cases}
					C\sum_{m=1}^{\infty}m^{-r(1/q-1/2)+\alpha-2}<\infty, & \text{ if $\alpha q\ge 2$}\\
					C\sum_{m=1}^{\infty}m^{-r(q\alpha/2+1/q-3/2)+\alpha-2}<\infty, &\text{if $\alpha q<2$,}
				\end{cases}
			\end{eqnarray*}
			\begin{eqnarray*}
				\sum_{m=1}^{\infty}m^{1/q-1}m^{-\min\{2,q\alpha\}r/2}\sum_{n=1}^{m}n^{r/2+\alpha-2-1/q}\le \begin{cases}
					C\sum_{m=1}^{\infty}m^{-r/2+\alpha-2}<\infty, &\text{ if $q\alpha\ge 2$,}\\
					C\sum_{m=1}^{\infty}m^{-r(q\alpha-1)/2+\alpha-2}<\infty, &\text{ if $q\alpha<2$,}
				\end{cases}
			\end{eqnarray*}
			and hence we obtain
			$$
			\sum_{n=1}^{\infty}n^{\alpha-2}\vv\left\{\max_{1\le k\le n}\left|\sum_{j=1}^{k}X_j\right|>\varepsilon n^{1/q}\right\}<\infty \mbox{ for any $\varepsilon>0$,}
			$$
			which complements that of Theorems 3.1 of Chen and Wu \cite{chen2022complete}.
		\end{rmk}
	\begin{rmk}\label{rmk03}
		The subtle difference between the classic probability space and the sub-linear expectation space here could be implied heuristically by the fact that although $\ee[\alpha Y_1+\beta Y_2]=\alpha \ee[Y_1]+\beta\ee[Y_2]$ holds for $\ee(-Y)=\ee(Y)=0$, but $\ee[\alpha Y_1^r+\beta Y_2^t]\neq\alpha\ee[Y_1^r]+\beta\ee[Y_2^s]$ for any $r\neq1$, $s\neq 1$, $\alpha,\beta\in \rr$.
	\end{rmk}
	
	As in Zhang and Ding \cite{zhang2017further}, assumptions \rm(C1)-\rm(C6) could be fulfilled by many sequences, the interested reader could refer to  Zhang and Ding \cite{zhang2017further} for more examples of similar sequences.

	\begin{thm}\label{thm2.2}Suppose $l$ is a function slowly varying at infinity. Assume that $\{X_n=\sum_{i=-\infty}^{\infty}a_iY_{i+n},n\ge 1\}$ is a moving average process generated by a sequence of negatively dependent random variables $\{Y_i,-\infty<i<\infty\}$ with $\ee(Y_i)=\ee(-Y_i)=0$, where $\{a_i,-\infty<i<\infty\}$ is an absolutely summable sequence of real numbers. Suppose that $\{g(n),n\ge 1\}$ and $\{f(n),n\ge 1\}$ are two sequences of positive constants with $f(n)\uparrow \infty$, and $\{\Psi_n(t),n\ge 1\}$ is a sequence of even and nonnegative functions such that, for each $n\ge 1$, $t>0$, $\Psi_n(t)>0$. Suppose that
		\begin{equation}\label{2.3}
			\frac{\Psi_n(t)}{|t|^p}\uparrow, \frac{\Psi_n(t)}{|t|^q}\downarrow, \mbox{  as  $|t|\uparrow$,}
		\end{equation}
		for some $1\le p<q<2$, and
		\begin{equation}\label{2.4}
			\sum_{n=1}^{\infty}g(n)l(n)\sum_{i=j+1}^{j+n}\frac{C_{\vv}(\Psi_i(Y_i))}{\Psi_i(f(n))}<\infty, \sum_{i=j+1}^{j+n}\frac{C_{\vv}(\Psi_i(Y_i))}{\Psi_i(f(n))}\rightarrow 0, \mbox{  as $n\rightarrow \infty$,}
		\end{equation}
		for any $j\ge 0$. Then for all $\varepsilon>0$,
		\begin{equation}\label{2.5}
			\sum_{n=1}^{\infty}g(n)l(n)\vv\left\{\max_{1\le k\le n}\left|\sum_{j=1}^{k}X_j\right|>\varepsilon f(n)\right\}<\infty.
		\end{equation}
	\end{thm}
	\begin{rmk}\label{rmk00}
		In Theorem \ref{thm2.2}, let $1\le s<2$, $1<\alpha<2/s$, $f(n)=n^{1/s}$, $g(n)=n^{\alpha-2}$, $l(n)=1$, $p=s\alpha$, $1\le p<q< 2$, $C_{\vv}\{|Y|^{q}\}<\infty$, $\Psi_n(x)=|x|^{q}$. Then assumptions of Theorem \ref{thm2.2} holds, and hence we obtain
		$$
		\sum_{n=1}^{\infty}n^{\alpha-2}\vv\left\{\max_{1\le k\le n}\left|\sum_{j=1}^{k}X_j\right|>\varepsilon n^{1/s}\right\}<\infty \mbox{ for any $\varepsilon>0$,}
		$$
		which complements that of Theorems 3.1 of Chen and Wu \cite{chen2022complete}.
	\end{rmk}
	
	\section{Proof of major results}
	
	\begin{proof}[\emph{ Proof of Theorem \ref{thm2.1}}] Write $Y_{xj}=-xI\{Y_j<-x\}+Y_jI\{|Y_j|\le x\}+xI\{Y_j>x\}$, $Y_{xj}'=Y_j-Y_{xj}$,  $Y_{x}=-xI\{Y<-x\}+YI\{|Y|\le x\}+xI\{Y>x\}$, $Y_{x}'=Y-Y_{x}$. Note that $\sum_{k=1}^{n}X_k=\sum_{i=-\infty}^{\infty}a_i\sum_{j=i+1}^{i+n}Y_j$. Notice that $\sum_{i=-\infty}^{\infty}|a_i|<\infty$, $\ee(Y_i)=0$, $C_{\vv}\left(|Y|^p(1\vee l(f^{-1}(|Y|)))\right)<\infty$, then by assumption \rm(C1), $|\ee(X)-\ee(Y)|\le \ee|X-Y|$, and Lemma 4.5 of Zhang \cite{zhang2016exponential}, for any $x>f(n)$, we see that
		$$
		\begin{aligned}
			&x^{-1}\max_{1\le k\le n}\left|\sum_{i-\infty}^{\infty}a_i\sum_{j=i+1}^{i+k}\ee(Y_{xj})\right|=x^{-1}\max_{1\le k\le n}\left|\sum_{i-\infty}^{\infty}a_i\sum_{j=i+1}^{i+k}[\ee(Y_{xj})-\ee(Y_j)]\right| \\
			&\quad\le x^{-1}\max_{1\le k\le n}C\sum_{i=-\infty}^{\infty}|a_i|\sum_{j=i+1}^{i+k}\ee|Y_{xj}'|\le Cx^{-1}\sum_{i=-\infty}^{\infty}|a_i|\sum_{j=i+1}^{i+n}\ee|Y_{xj}'| \\
			&\quad =Cx^{-1}\sum_{i=-\infty}^{\infty}|a_i|\sum_{j=i+1}^{i+n}\ee|Y_{x}'|\\
			&\quad\le Cx^{-1}\sum_{i=-\infty}^{\infty}|a_i|\sum_{j=i+1}^{i+n}C_{\vv}\{|Y_{x}'|\}\quad\le Cx^{-1}\sum_{i=-\infty}^{\infty}|a_i|\sum_{j=i+1}^{i+n}C_{\vv}\left(|Y|I\{|Y|>x\}\right)\\
			&\quad \le Cx^{-1}nC_{\vv}\left(|Y|I\{|Y|>x\}\right) =Cx^{-1}n\int_{0}^{\infty}\vv\left(|Y|I\{|Y|>x\}>y\right)\dif y\\
			&\quad =Cx^{-1}n\left[\int_{0}^{x}\vv\left(|Y|>x\right)\dif y+\int_{x}^{\infty}\vv\left(|Y|>y\right)\dif y\right]\le Cx^{-1}n\Bigg[\vv\left(|Y|>x\right)x\\
			&\quad\quad +\left.\int_{x}^{\infty}1\cdot\vv\left(|Y|>y\right)\dif y\right] \le Cx^{-1}n\left[\int_{0}^{x}\frac{py^{p-1}}{x^{p-1}}\vv\left(|Y|>x\right)\dif y+\int_{x}^{\infty}\frac{py^{p-1}}{x^{p-1}}\vv\left(|Y|>y\right)\dif y\right]\\
			&\quad \le cnx^{-p}C_{\vv}\left(|Y|^pI\{|Y|>x\}\right)\le C\frac{n}{f^p(n)}C_{\vv}\left(|Y|^pI\{|Y|>x\}\right)\rightarrow 0,
			\mbox{ as $x\rightarrow\infty$.}
		\end{aligned}
		$$
		Hence, we obtain
		\[
		x^{-1}\max_{1\le k\le n}\left|\sum_{i=-\infty}^{\infty}a_i\sum_{j=i+1}^{i+k}\ee(Y_{xj})\right|<\varepsilon/4,
		\]
		for any $\varepsilon>0$ and $x>f(n)$ large sufficiently. Therefore we conclude that
		\begin{eqnarray}\label{4.1}
			\nonumber && \sum_{n=1}^{\infty}\frac{g(n)l(n)}{f(n)}C_{\vv}\left\{\left(\max_{1\le k\le n}\left|\sum_{j=1}^{k}X_j\right|-\varepsilon f(n)\right)^{+}\right\}
		\end{eqnarray}
		\begin{eqnarray}
			\nonumber&&\quad \le \sum_{n=1}^{\infty}\frac{g(n)l(n)}{f(n)}\int_{\varepsilon f(n)}^{\infty}\vv\left\{\max_{1\le k\le n}\left|\sum_{j=1}^{k}X_j\right|>x\right\}\dif x\\
			\nonumber&&\quad\le \sum_{n=1}^{\infty}\frac{g(n)l(n)}{f(n)}\int_{ f(n)}^{\infty}\vv\left\{\max_{1\le k\le n}\left|\sum_{j=1}^{k}X_j\right|>\varepsilon x\right\}\dif x\\
			\nonumber&&\quad \le \sum_{n=1}^{\infty}\frac{g(n)l(n)}{f(n)}\int_{ f(n)}^{\infty}\vv\left\{\max_{1\le k\le n}\left|\sum_{i=-\infty}^{\infty}a_i\sum_{j=i+1}^{i+k}(Y_j-Y_{xj})\right|>\varepsilon x/2\right\}\dif x\\
			\nonumber &&\quad \quad +\sum_{n=1}^{\infty}\frac{g(n)l(n)}{f(n)}\int_{ f(n)}^{\infty}\vv\left\{\max_{1\le k\le n}\left|\sum_{i=-\infty}^{\infty}a_i\sum_{j=i+1}^{i+k}(Y_{xj}-\ee Y_{xj})\right|>\varepsilon x/4\right\}\dif x\\
			&&\quad =: \Rmnum{1}_1+\Rmnum{1}_2.
		\end{eqnarray}
		
		Now we want to establish $\Rmnum{1}_1<\infty$. Obviously $|Y_j-Y_{xj}|\le |Y_j|I\{|Y_j|>x\}$, write $H(m)=\sum_{k=1}^{m}\log\frac{f(k+1)}{f(k)}\sum_{n=1}^{k}\frac{ng(n)l(n)}{f(n)}$, $H(0)=0$ below, by Markov's inequality under sub-linear expectations, assumptions \rm(C1) and \rm(C2), and the proof of Lemma 2.2 in Zhong and Wu \cite{zhong2017complete}, we see that
		\begin{eqnarray*}
			\nonumber \Rmnum{1}_1&\le&C  \sum_{n=1}^{\infty}\frac{g(n)l(n)}{f(n)}\int_{ f(n)}^{\infty}x^{-1}\ee^{*}\max_{1\le k\le n}\left|\sum_{i=-\infty}^{\infty}a_i\sum_{j=i+1}^{i+k}(Y_j-Y_{xj})\right|\dif x\\
			\nonumber &\le& C \sum_{n=1}^{\infty}\frac{g(n)l(n)}{f(n)}\int_{ f(n)}^{\infty}x^{-1}\sum_{i=-\infty}^{\infty}|a_i|\sum_{j=i+1}^{i+n}\ee^{*}\left|Y_j-Y_{xj}\right|\dif x\\
			\nonumber &\le& C \sum_{n=1}^{\infty}\frac{g(n)l(n)}{f(n)}\int_{ f(n)}^{\infty}x^{-1}\sum_{i=-\infty}^{\infty}|a_i|\sum_{j=i+1}^{i+n}\ee\left|Y_j-Y_{xj}\right|\dif x\\
			\nonumber&=& C \sum_{n=1}^{\infty}\frac{ng(n)l(n)}{f(n)}\int_{ f(n)}^{\infty}x^{-1}\sum_{i=-\infty}^{\infty}|a_i|\ee[|Y_x'|]\dif x\\
			\nonumber &\le &C\sum_{n=1}^{\infty}\frac{ng(n)l(n)}{f(n)}\int_{ f(n)}^{\infty}x^{-1}C_{\vv}\left(|Y|I\{|Y|>x\}\right)\dif x\\
			\nonumber &=&C\sum_{n=1}^{\infty}\frac{ng(n)l(n)}{f(n)}\sum_{m=n}^{\infty}\int_{f(m)}^{f(m+1)}x^{-1}C_{\vv}\left(|Y|I\{|Y|>x\}\right)\dif x\\
			\nonumber &\le&C\sum_{n=1}^{\infty}\frac{ng(n)l(n)}{f(n)}\sum_{m=n}^{\infty}\log\frac{f(m+1)}{f(m)}C_{\vv}\left(|Y|I\{|Y|>f(m)\}\right)\\
			\nonumber &=&C\sum_{m=1}^{\infty}\log\frac{f(m+1)}{f(m)}C_{\vv}\left(|Y|I\{|Y|>f(m)\}\right)\sum_{n=1}^{m}\frac{ng(n)l(n)}{f(n)}\\
			\nonumber &=&C\sum_{m=1}^{\infty}\log\frac{f(m+1)}{f(m)}\sum_{n=1}^{m}\frac{ng(n)l(n)}{f(n)}\left[\int_{0}^{f(m)}\vv\left(|Y|>f(m)\right)\dif y+\int_{f(m)}^{\infty}\vv(|Y|>y)\dif y\right]\\
			\nonumber &=&C\sum_{m=1}^{\infty}\log\frac{f(m+1)}{f(m)}\sum_{n=1}^{m}\frac{ng(n)l(n)}{f(n)}f(m)\vv\left(|Y|>f(m)\right)
		\end{eqnarray*}
		\begin{eqnarray}
			\nonumber&& +C\sum_{m=1}^{\infty}\log\frac{f(m+1)}{f(m)}\sum_{n=1}^{m}\frac{ng(n)l(n)}{f(n)}\int_{f(m)}^{\infty}\vv(|Y|>y)\dif y\\
			\nonumber &=&C\lim_{k\rightarrow \infty}\sum_{m=1}^{k}\left(H(m)-H(m-1)\right)f(m)\vv\left(|Y|>f(m)\right)\\
			\nonumber &&+C\sum_{m=1}^{\infty}\log\frac{f(m+1)}{f(m)}\sum_{n=1}^{m}\frac{ng(n)l(n)}{f(n)}\sum_{k=m}^{\infty}\int_{f(k)}^{f(k+1)}\vv(|Y|>y)\dif y\\
			\nonumber &\le&C\lim_{k\rightarrow\infty}\left(\sum_{m=1}^{k-1}H(m)\left[f(m)\vv\left(|Y|>f(m)\right)-f(m+1)\vv\left(|Y|>f(m+1)\right)\right]\right.\\
			\nonumber&&+H(k)f(k)\vv\left(|Y|>f(k)\right)\Bigg)+C\sum_{k=1}^{\infty}\sum_{m=1}^{k}\log\frac{f(m+1)}{f(m)}\sum_{n=1}^{m}\frac{ng(n)l(n)}{f(n)}\int_{f(k)}^{f(k+1)}\vv(|Y|>y)\dif y\\
			\nonumber &\le &C\lim_{k\rightarrow\infty}\left(C\sum_{m=1}^{k-1}f^{p-1}(m)l(m)\left[f(m)\vv\left(|Y|>f(m)\right)-f(m+1)\vv\left(|Y|>f(m+1)\right)\right]\right.\\ \nonumber&&+Cf^{p}(k)l(k)\vv\left(|Y|>f(k)\right)\Bigg)+C\sum_{k=1}^{\infty}f^{p-1}(k)l(k)\int_{f(k)}^{f(k+1)}\vv(|Y|>y)\dif y\\
			\nonumber &\le& C\lim_{k\rightarrow\infty}\sum_{m=1}^{k}\left(f^{p-1}(m)l(m)-f^{p-1}(m-1)l(m-1)\right)f(m)\vv\left(|Y|^pl(f^{-1}(|Y|))>f^p(m)l(m)\right)\\
			\nonumber &&+C\int_{f(1)}^{\infty}\vv\left(|Y|>y\right)y^{p-1}l(f^{-1}(y))\dif y\\
			&\le &C C_{\vv}\left(|Y|^p\left(1\vee l(f^{-1}(|Y|))\right)\right)<\infty.
		\end{eqnarray}
		Therefore it remains to prove that $\Rmnum{1}_2<\infty$. By Markov's inequality under sub-linear expectations, H\"{o}lder's inequality, and Lemma \ref{lem01}, for $r>\max\{2,p\}$, we see that
		\begin{eqnarray}\label{4.3}
			\nonumber \Rmnum{1}_2&\le&C  \sum_{n=1}^{\infty}\frac{g(n)l(n)}{f(n)}\int_{ f(n)}^{\infty}x^{-r}\ee^{*}\max_{1\le k\le n}\left|\sum_{i=-\infty}^{\infty}a_i\sum_{j=i+1}^{i+k}(Y_{xj}-\ee Y_{xj})\right|^{r}\dif x\\
			\nonumber &\le& C \sum_{n=1}^{\infty}\frac{g(n)l(n)}{f(n)}\int_{ f(n)}^{\infty}x^{-r}\ee^{*}\left[\sum_{i=-\infty}^{\infty}\left(|a_i|^{\frac{r-1}{r}}\right)\left(|a_i|^{1/r}\max_{1\le k\le n}\left|\sum_{j=i+1}^{i+k}(Y_{xj}-\ee Y_{xj})\right|\right)\right]^r\dif x\\
			\nonumber &\le& C\sum_{n=1}^{\infty}\frac{g(n)l(n)}{f(n)}\int_{ f(n)}^{\infty}x^{-r}\left(\sum_{i=-\infty}^{\infty}|a_i|\right)^{r-1}\left(\sum_{i=-\infty}^{\infty}|a_i|\ee^{*}\max_{1\le k\le n}\left|\sum_{j=i+1}^{i+k}(Y_{xj}-\ee Y_{xj})\right|^r\right)\dif x\\
			\nonumber &=& C\sum_{n=1}^{\infty}\frac{g(n)l(n)}{f(n)}\int_{ f(n)}^{\infty}x^{-r}\left(\sum_{i=-\infty}^{\infty}|a_i|\right)^{r-1}\left(\sum_{i=-\infty}^{\infty}|a_i|\ee\max_{1\le k\le n}\left|\sum_{j=i+1}^{i+k}(Y_{xj}-\ee Y_{xj})\right|^r\right)\dif x\\
			\nonumber &\le&C\sum_{n=1}^{\infty}\frac{g(n)l(n)}{f(n)}\int_{ f(n)}^{\infty}x^{-r}\sum_{i=-\infty}^{\infty}|a_i|\sum_{j=i+1}^{i+n}\ee\left|Y_{xj}-\ee Y_{xj}\right|^r\dif x
		\end{eqnarray}
		\begin{eqnarray}\label{4.3}
			\nonumber &&+C\sum_{n=1}^{\infty}\frac{g(n)l(n)}{f(n)}\int_{ f(n)}^{\infty}x^{-r}\sum_{i=-\infty}^{\infty}|a_i|\left(\sum_{j=i+1}^{i+n}\ee\left|Y_{xj}-\ee Y_{xj}\right|^2\right)^{r/2}\dif x\\
			\nonumber	&&+C\sum_{n=1}^{\infty}\frac{g(n)l(n)}{f(n)}\int_{ f(n)}^{\infty}x^{-r}\sum_{i=-\infty}^{\infty}|a_i|\left(\sum_{j=i+1}^{i+n}\left(-\ee(-Y_{xj})-\ee(Y_{xj})\right)^{-}\right)^{r}\dif x\\
			&=&:\Rmnum{1}_{21}+\Rmnum{1}_{22}+\Rmnum{1}_{23}.
		\end{eqnarray}
		For $\Rmnum{1}_{21}$, by $C_r$ inequality, assumptions \rm(C1), \rm(C3), \rm(C4), and Lemma 4.5 of Zhang \cite{zhang2016exponential}, write 
		$$\hat{H}(m)=\sum_{k=1}^{m}\left(f(k+1)-f(k)\right)\sum_{n=1}^{k}\frac{ng(n)l(n)}{f(n)},m\ge 1,$$
		$$
		\hat{H}(0)=0,
		$$
		we conclude that
		\begin{eqnarray}\label{4.4}
			\nonumber \Rmnum{1}_{21}&\le&C  \sum_{n=1}^{\infty}\frac{ng(n)l(n)}{f(n)}\int_{ f(n)}^{\infty}x^{-r}\sum_{i=-\infty}^{\infty}|a_i|\sum_{j=i+1}^{i+n}\ee|Y_{xj}|^r\dif x= C  \sum_{n=1}^{\infty}\frac{g(n)l(n)}{f(n)}\int_{ f(n)}^{\infty}x^{-r}\ee|Y_{x}|^r\dif x\\
			\nonumber &\le&C\sum_{n=1}^{\infty}\frac{ng(n)l(n)}{f(n)}\int_{ f(n)}^{\infty}x^{-r}C_{\vv}\left(|Y_{x}|^r\right)\dif x\\
			\nonumber &\le&C\sum_{n=1}^{\infty}\frac{ng(n)l(n)}{f(n)}\sum_{m=n}^{\infty}\int_{f(m)}^{f(m+1)}\left[x^{-r}C_{\vv}\left(|Y|^rI\{|Y|\le x\}\right)+\vv\left(|Y|> x\right)\right]\dif x\\
			\nonumber &\le&C\sum_{m=1}^{\infty}\left(f^{1-r}(m)-f^{1-r}(m+1)\right)\sum_{n=1}^{m}\frac{ng(n)l(n)}{f(n)}\sum_{k=1}^{m}\int_{f(k)}^{f(k+1)}\vv\left(|Y|>x\right)x^{r-1}r\dif x\\
			\nonumber &&+C\sum_{m=1}^{\infty}\left[f(m+1)-f(m)\right]\vv\left(|Y|>f(m)\right)\sum_{n=1}^{m}\frac{ng(n)l(n)}{f(n)}\\
			\nonumber &\le&C\sum_{k=1}^{\infty}\int_{f(k)}^{f(k+1)}\vv\left(|Y|>x\right)x^{r-1}r\dif x\sum_{m=k}^{\infty}\left(f^{1-r}(m)-f^{1-r}(m+1)\right)\sum_{n=1}^{m}\frac{ng(n)l(n)}{f(n)}\\
			\nonumber &&+C \sum_{m=1}^{\infty}\left(\hat{H}(m)-\hat{H}(m-1)\right)\vv\left(|Y|>f(m)\right)\\
			\nonumber &\le&C\sum_{k=1}^{\infty}f^{p-r}(k)l(k)\int_{f(k)}^{f(k+1)}\vv\left(|Y|>x\right)x^{r-1}r\dif x\\
			\nonumber &&+\lim_{k\rightarrow \infty}C\sum_{m=1}^{k-1}\hat{H}(m)\left(\vv\left(|Y|>f(m)\right)-\vv\left(|Y|>f(m+1)\right)\right)+\hat{H}(k)\vv\left(|Y|>f(k)\right)\\
			\nonumber &\le&C\sum_{k=1}^{\infty}\int_{f(k)}^{f(k+1)}\vv\left(|Y|>x\right)x^{p-1}l(f^{-1}(x))\dif x\\
			\nonumber &&+\lim_{k\rightarrow \infty}C\sum_{m=1}^{k-1}f^p(m)l(m)\left(\vv\left(|Y|>f(m)\right)-\vv\left(|Y|>f(m+1)\right)\right)+f^p(k)l(k)\vv\left(|Y|>f(k)\right)\\
			\nonumber &\le&C C_{\vv}\left(|Y|^p\left(1\vee l(f^{-1}(|Y|))\right)\right)
		\end{eqnarray}
		\begin{eqnarray}\label{4.4}	
			\nonumber &&+\lim_{k\rightarrow \infty}C\sum_{m=1}^{k}\left(f^p(m)l(m)-f^p(m-1)l(m-1)\right)\vv\left(|Y|>f(m)\right)\\
			\nonumber
			&\le&C C_{\vv}\left(|Y|^p\left(1\vee l(f^{-1}(|Y|))\right)\right)<\infty.
		\end{eqnarray}
		
		Next, we want to establish $\Rmnum{1}_{22}<\infty$. By $C_r$ ienquality, Lemma 4.5 of Zhang \cite{zhang2016exponential}, and assumptions \rm(C1), \rm(C5), and \rm(C6), we see that
		\begin{eqnarray}\label{4.5}
			\nonumber \Rmnum{1}_{22}&\le&C  \sum_{n=1}^{\infty}\frac{n^{r/2}g(n)l(n)}{f(n)}\int_{ f(n)}^{\infty}x^{-r}\left(\ee|Y_{x1}|^2\right)^{r/2}\dif x= C  \sum_{n=1}^{\infty}\frac{n^{r/2}g(n)l(n)}{f(n)}\int_{ f(n)}^{\infty}x^{-r}\left(\ee|Y_{x}|^2\right)^{r/2}\dif x\\ \nonumber&\le&C\sum_{n=1}^{\infty}\frac{n^{r/2}g(n)l(n)}{f(n)}\int_{ f(n)}^{\infty}x^{-r}\left(C_{\vv}\left\{|Y_{x}|^2\right\}\right)^{r/2}\dif x\\
			\nonumber&=&C\sum_{n=1}^{\infty}\frac{n^{r/2}g(n)l(n)}{f(n)}\int_{ f(n)}^{\infty}x^{-r}\left[\int_{0}^{\infty}\vv\left(|Y|^2I\{|Y|\le x\}+x^2I\{|Y|>x\}>y\right)\dif y\right]^{r/2}\dif x\\
			\nonumber&\le&C\sum_{n=1}^{\infty}\frac{n^{r/2}g(n)l(n)}{f(n)}\int_{ f(n)}^{\infty}x^{-r}\left[\int_{0}^{\infty}\left[\vv\left(|Y|^2I\{|Y|\le x\}>y/2\right)+\vv\left(x^2I\{|Y|>x\}>y/2\right)\right]\dif y\right]^{r/2}\dif x \\
			\nonumber&\le&C\sum_{n=1}^{\infty}\frac{n^{r/2}g(n)l(n)}{f(n)}\int_{ f(n)}^{\infty}x^{-r}\left[C_{\vv}\{|Y|^2I\{|Y|\le x\}\}+x^2\vv\left(|Y|>x\right)\right]^{r/2}\dif x \\
			\nonumber &\le&C\sum_{n=1}^{\infty}\frac{n^{r/2}g(n)l(n)}{f(n)}\sum_{m=n}^{\infty}\int_{f(m)}^{f(m+1)}\left[x^{-r}\left(C_{\vv}\{|Y|^2I\{|Y|\le x\}\}\right)^{r/2}+\vv^{r/2}\left(|Y|>x\right)\right]\dif x\\
			\nonumber &\le&C\sum_{n=1}^{\infty}\frac{n^{r/2}g(n)l(n)}{f(n)}\sum_{m=n}^{\infty}\left[f^{1-r}(m)-f^{1-r}(m+1)\right]\left(C_{\vv}\{|Y|^2I\{|Y|\le f(m+1)\}\}\right)^{r/2}\\
			\nonumber &&+C\sum_{n=1}^{\infty}\frac{n^{r/2}g(n)l(n)}{f(n)}\sum_{m=n}^{\infty}\left[f(m+1)-f(m)\right]\vv^{r/2}\left(|Y|>f(m)\right)\\
			\nonumber&=&C\sum_{m=1}^{\infty}\left[f^{1-r}(m)-f^{1-r}(m+1)\right]\left(C_{\vv}\{|Y|^2I\{|Y|\le f(m+1)\}\}\right)^{r/2}\sum_{n=1}^{m}\frac{n^{r/2}g(n)l(n)}{f(n)}\\
			\nonumber &&+\sum_{m=1}^{\infty}\left[f(m+1)-f(m)\right]\vv^{r/2}\left(|Y|>f(m)\right)\sum_{n=1}^{m}\frac{n^{r/2}g(n)l(n)}{f(n)}\\
			\nonumber&\le&C\sum_{m=1}^{\infty}\left[f^{1-r}(m)-f^{1-r}(m+1)\right]\sum_{n=1}^{m}\frac{n^{r/2}g(n)l(n)}{f(n)}f^{\max\{0,2-p\}r/2}(m+1)\left(C_{\vv}\left\{|Y|^{\min\{p,2\}}\right\}\right)^{r/2}\\
			&&+C\sum_{m=1}^{\infty}\left[f(m+1)-f(m)\right]\sum_{n=1}^{m}\frac{n^{r/2}g(n)l(n)}{f(n)}f^{-\min\{2,p\}r/2}(m)\left(\ee|Y|^{\min\{2,p\}}\right)^{r/2}<\infty.
		\end{eqnarray}
		Finally, we aim to prove $\Rmnum{1}_{23}<\infty$. By $\ee(Y)=\ee(-Y)=0$, $|\ee(X)-\ee(Y)|\le \ee|X-Y|$, \rm(C1), \rm(C5), and Lemma 4.5 of Zhang \cite{zhang2016exponential}, we see that
		\begin{eqnarray}\label{4.6}
			\nonumber	\Rmnum{1}_{23}&\le& C\sum_{n=1}^{\infty}\frac{n^rg(n)l(n)}{f(n)}\int_{ f(n)}^{\infty}x^{-r}\left(\ee(-Y_{x1})+\ee(Y_{x1})\right)^{r}\dif x\\
			\nonumber&=&C\sum_{n=1}^{\infty}\frac{n^rg(n)l(n)}{f(n)}\int_{ f(n)}^{\infty}x^{-r}\left(\ee(-Y_{x})+\ee(Y_{x})\right)^{r}\dif x
		\end{eqnarray}
		\begin{eqnarray}\label{4.6}
			\nonumber&\le&C\sum_{n=1}^{\infty}\frac{n^rg(n)l(n)}{f(n)}\int_{ f(n)}^{\infty}x^{-r}\left(\ee|Y_x|\right)^{r}\dif x\le \sum_{n=1}^{\infty}\frac{n^rg(n)l(n)}{f(n)}\int_{ f(n)}^{\infty}x^{-r}\left(C_{\vv}\left\{|Y_x|\right\}\right)^{r}\dif x\\
			\nonumber	&\le& C\sum_{n=1}^{\infty}\frac{n^rg(n)l(n)}{f(n)}\sum_{m=n}^{\infty}\int_{f(m)}^{f(m+1)}x^{-r}\left(C_{\vv}[|Y|I\{|Y|>f(m)\}]\right)^r\dif x\\
			\nonumber	&\le& C\sum_{m=1}^{\infty}\left(f^{1-r}(m)-f^{1-r}(m+1)\right)\sum_{n=1}^{m}\frac{n^rg(n)l(n)}{f(n)}\left(\frac{C_{\vv}\{|Y|^pI\{|Y|>f(m)\}\}}{f^{p-1}(m)}\right)^r\\
			\nonumber	&\le& C\sum_{m=1}^{\infty}\left(f^{1-r}(m)-f^{1-r}(m+1)\right)f^{(2-p)r/2}(m)\sum_{n=1}^{m}\frac{n^{r/2}g(n)l(n)}{f(n)}\frac{n^{r/2}}{f^{pr/2}(m)}\left(C_{\vv}\left(|Y|^p\right)\right)^r\\
			\nonumber	&\le& C\sum_{m=1}^{\infty}\left(f^{1-r}(m)-f^{1-r}(m+1)\right)f^{(2-p)r/2}(m)\sum_{n=1}^{m}\frac{n^{r/2}g(n)l(n)}{f(n)}\\
			&\le& C\sum_{m=1}^{\infty}\left(f^{1-r}(m)-f^{1-r}(m+1)\right)f^{\max\{0,(2-p)r/2\}}(m+1)\sum_{n=1}^{m}\frac{n^{r/2}g(n)l(n)}{f(n)}<\infty.
		\end{eqnarray}
		Therefore, combining (\ref{4.1})-(\ref{4.6}) results in (\ref{2.1}). The proof is complete.
	\end{proof}
	
	\begin{proof}[\emph{ Proof of Theorem \ref{thm2.2}}] Write $Y_{nj}=-f(n)I\{Y_j<f(n)\}+Y_jI\{|Y_j|\le f(n)\}+f(n)I\{Y_j>f(n)\}$, $Y_{nj}'=Y_{nj}-Y_j$. Obviously $\sum_{j=1}^{k}X_j=\sum_{i=-\infty}^{\infty}a_i\sum_{j=i+1}^{i+k}Y_j$. Observing that $\sum_{i=-\infty}^{\infty}|a_i|<\infty$ and $\ee(Y_j)=0$, $|\ee(X)-\ee(Y)|\le \ee|X-Y|$, then by $C_{\vv}\{\alpha |X|\}=\int_{0}^{\infty}\lambda\vv\{|X|>x/\lambda\}\dif x/\lambda=\lambda C_{\vv}\{|X|\}$, $\forall \lambda>0$, (\ref{2.3}) and (\ref{2.4}), we see that
		\begin{eqnarray*}
			&&\frac{1}{f(n)}\max_{1\le k\le n}\left|\sum_{i=-\infty}^{\infty}a_i\sum_{j=i+1}^{i+k}\ee(Y_{nj})\right|\le \frac{1}{f(n)}\sum_{i=-\infty}^{\infty}|a_i|\sum_{j=i+1}^{i+n}|\ee(Y_{nj})-\ee(Y_j)|\\
			&&\quad\le \frac{1}{f(n)}\sum_{i=-\infty}^{\infty}|a_i|\sum_{j=i+1}^{i+n}C_{\vv}\{|Y_{nj}'|\}\le\frac{1}{f(n)}\sum_{i=-\infty}^{\infty}|a_i|\sum_{j=i+1}^{i+n}C_{\vv}\{|Y_j|I\{|Y_j|>f(n)\}\}  \\
			&&\quad=  \sum_{i=-\infty}^{\infty}|a_i|\sum_{j=i+1}^{i+n}C_{\vv}\left\{\frac{|Y_j|I\{|Y_j|>f(n)\}}{f(n)}\right\}\le \sum_{i=-\infty}^{\infty}|a_i|\sum_{j=i+1}^{i+n}C_{\vv}\left\{\frac{|Y_j|^pI\{|Y_j|>f(n)\}}{f^p(n)}\right\}\\
			&&\quad\le C\sum_{i=-\infty}^{\infty}|a_i|\sum_{j=i+1}^{i+n}C_{\vv}\left[\frac{\Psi_j(Y_j)}{\Psi_j(f(n))}\left[I\left\{|Y_j|\le f(n)\right\}+I\left\{|Y_j|>f(n)\right\}\right]\right]\rightarrow 0, \mbox{  as $n\rightarrow\infty$.}
		\end{eqnarray*}
		Therefore for $n$ sufficiently large and any $\varepsilon>0$, we have
		\[
		\frac{1}{f(n)}\max_{1\le k\le n}\left|\sum_{i=-\infty}^{\infty}a_i\sum_{j=i+1}^{i+k}\ee(Y_{nj})\right|<\frac{\varepsilon}{4}.
		\]
		Then we conclude that
		\begin{eqnarray*}
			&&\sum_{n=1}^{\infty}g(n)l(n)\vv\left\{\max_{1\le k\le n}\left|\sum_{j=1}^{k}X_j\right|>\varepsilon f(n)\right\}\\
			&&\quad \le C\sum_{n=1}^{\infty}g(n)l(n)\vv\left\{\max_{1\le k\le n}\left|\sum_{i=-\infty}^{\infty}a_i\sum_{j=i+1}^{i+k}\left(Y_j-Y_{nj}\right)\right|>\varepsilon f(n)/2\right\}
		\end{eqnarray*}
		\begin{eqnarray*}
			&&\quad\quad +C\sum_{n=1}^{\infty}g(n)l(n)\vv\left\{\max_{1\le k\le n}\left|\sum_{i=-\infty}^{\infty}a_i\sum_{j=i+1}^{i+k}\left(Y_{nj}-\ee Y_{nj}\right)\right|>\varepsilon f(n)/4\right\}\\
			&&\quad=:J_1+J_2.
		\end{eqnarray*}
		By Markov's inequality under sub-linear expectations, Lemma 4.5 of Zhang \cite{zhang2016exponential}, (\ref{2.3}), and (\ref{2.4}),  we see that
		$$
		\begin{aligned}
			J_1&\le C\sum_{n=1}^{\infty}\frac{g(n)l(n)}{f(n)}\ee^{*}\max_{1\le k\le n}\left|\sum_{i=-\infty}^{\infty}a_i\sum_{j=i+1}^{i+k}\left(Y_j-Y_{nj}\right)\right|\le C\sum_{n=1}^{\infty}\frac{g(n)l(n)}{f(n)}\sum_{i=-\infty}^{\infty}|a_i|\sum_{j=i+1}^{i+n}\ee|Y_{nj}'|\\
			&\le C\sum_{n=1}^{\infty}\frac{g(n)l(n)}{f(n)}\sum_{i=-\infty}^{\infty}|a_i|\sum_{j=i+1}^{i+n}C_{\vv}\{|Y_{nj}'|\}\le C\sum_{n=1}^{\infty}\frac{g(n)l(n)}{f(n)}\sum_{i=-\infty}^{\infty}|a_i|\sum_{j=i+1}^{i+n}C_{\vv}\left\{|Y_j|I\{|Y_j|>f(n)\}\right\}\\
			&\le C\sum_{n=1}^{\infty}g(n)l(n)\sum_{i=-\infty}^{\infty}|a_i|\sum_{j=i+1}^{i+n}C_{\vv}\left\{\frac{|Y_j|I\{|Y_j|>f(n)\}}{f(n)}\right\}\le C\sum_{n=1}^{\infty}g(n)l(n)\sum_{i=-\infty}^{\infty}|a_i|\\
			&\quad\times\sum_{j=i+1}^{i+n}C_{\vv}\left\{\frac{|Y_j|^pI\{|Y_j|>f(n)\}}{f^p(n)}\right\}\le C\sum_{i=-\infty}^{\infty}|a_i|\sum_{n=1}^{\infty}g(n)l(n)\sum_{j=i+1}^{i+n}C_{\vv}\left(\frac{\Psi_j(Y_j)}{\Psi_j(f(n))}\right)<\infty.
		\end{aligned}
		$$
		By Markov's inequality under sub-linear expectations, H\"{o}lder's inequality, Lemma \ref{lem01}, $\ee(Y_i)=\ee(-Y_i)=0, i=1, 2, \ldots$, Proposition 1.3.7 of Peng \cite{peng2019nonlinear}, (\ref{2.3}), and (\ref{2.4}), and the proof of $J_1$, we obtain
		$$
		\begin{aligned}
			J_2&\le C\sum_{n=1}^{\infty}\frac{g(n)l(n)}{f^2(n)}\ee^{*}\max_{1\le k\le n}\left|\sum_{i=-\infty}^{\infty}a_i\sum_{j=i+1}^{i+k}\left(Y_{nj}-\ee Y_{nj}\right)\right|^2\\
			&\le C\sum_{n=1}^{\infty}\frac{g(n)l(n)}{f^2(n)}\sum_{i=-\infty}^{\infty}|a_i|\left(\sum_{j=i+1}^{i+n}\ee|Y_{nj}-\ee Y_{nj}|^2+\left(\sum_{j=i+1}^{i+n}\left[\ee(-Y_{nj})+\ee(Y_{nj})\right]\right)^2\right)\\
			&\le C\sum_{n=1}^{\infty}\frac{g(n)l(n)}{f^2(n)}\sum_{i=-\infty}^{\infty}|a_i|\sum_{j=i+1}^{i+n}\ee|Y_{nj}|^2\\
			&\quad +C\sum_{n=1}^{\infty}\frac{g(n)l(n)}{f^2(n)}\sum_{i=-\infty}^{\infty}|a_i|\left(\sum_{j=i+1}^{i+n}\left[|\ee(-Y_{nj}+Y_j)|+|\ee(Y_{nj}-Y_j)|\right]\right)^2\\
			&\le C\sum_{n=1}^{\infty}\frac{g(n)l(n)}{f^q(n)}\sum_{i=-\infty}^{\infty}|a_i|\sum_{j=i+1}^{i+n}\ee|Y_{nj}|^q +C\sum_{n=1}^{\infty}\frac{g(n)l(n)}{f^2(n)}\sum_{i=-\infty}^{\infty}|a_i|\left(\sum_{j=i+1}^{i+n}\ee|Y_{nj}'|\right)^2\\
			&\le C\sum_{i=-\infty}^{\infty}|a_i|\sum_{n=1}^{\infty}g(n)l(n)\sum_{j=i+1}^{i+n}\ee\frac{\Psi_j(Y_j)}{\Psi_j(f(n))} +C\sum_{n=1}^{\infty}g(n)l(n)\sum_{i=-\infty}^{\infty}|a_i|\left(\sum_{j=i+1}^{i+n}\ee\left[\frac{|Y_{nj}'|}{f(n)}\right]\right)^2\\
			&\le C+C\sum_{n=1}^{\infty}g(n)l(n)\sum_{i=-\infty}^{\infty}|a_i|\left(\sum_{j=i+1}^{i+n}C_{\vv}\left[\frac{|Y_{nj}'|}{f(n)}\right]\right)^2
		\end{aligned}
		$$
		$$
		\begin{aligned}
			&\le C+C\sum_{n=1}^{\infty}g(n)l(n)\sum_{i=-\infty}^{\infty}|a_i|\left(\sum_{j=i+1}^{i+n}C_{\vv}\left[\frac{|Y_j|I\{|Y_j|>f(n)\}}{f(n)}\right]\right)^2\\
			&\le C+C\sum_{n=1}^{\infty}g(n)l(n)\sum_{i=-\infty}^{\infty}|a_i|\left(\sum_{j=i+1}^{i+n}C_{\vv}\left[\frac{|Y_j|^pI\{|Y_j|>f(n)\}}{f^p(n)}\right]\right)^2\\
			&\le C+C\sum_{n=1}^{\infty}g(n)l(n)\sum_{i=-\infty}^{\infty}|a_i|\left(\sum_{j=i+1}^{i+n}C_{\vv}\left[\frac{\Psi_j(|Y_j|)}{\Psi_j(f(n))}\right]\right)^2\\
			&\le C+C\sum_{n=1}^{\infty}g(n)l(n)\sum_{i=-\infty}^{\infty}|a_i|\left(\sum_{j=i+1}^{i+n}\frac{C_{\vv}\{\Psi_j(|Y_j|)\}}{\Psi_j(f(n))}\right)^2\\
			&\le C+C\sum_{n=1}^{\infty}g(n)l(n)\sum_{j=i+1}^{i+n}\frac{C_{\vv}\{\Psi_j(|Y_j|)\}}{\Psi_j(f(n))}<\infty.
		\end{aligned}
		$$
		Hence the proof of Theorem \ref{thm2.2} is complete.
	\end{proof}
	


	
	
	
	
	

\end{document}